\documentclass[11pt]{article}
\usepackage{amsmath,amssymb}

\newtheorem{definition}{Definition}
\newtheorem{lemma}[definition]{Lemma}
\newtheorem{theorem}[definition]{Theorem}
\newtheorem{corollary}[definition]{Corollary}

\DeclareMathOperator{\ind}{ind}

\begin{document}

\title{Cusps and the family hyperbolic metric}         
\author{Scott A. Wolpert}        
\date{March 13, 2007}          

\maketitle

\begin{abstract}
The hyperbolic metric for the punctured unit disc in the Euclidean plane is {\em singular} at the origin.  A renormalization of the metric at the origin is provided by the Euclidean metric.  For Riemann surfaces there is a unique germ for the isometry class of a complete hyperbolic metric at a cusp.  The renormalization of the metric for the punctured unit disc provides a renormalization for a hyperbolic metric at a cusp.  For a holomorphic family of punctured Riemann surfaces the family of (co)tangent spaces along a puncture defines a {\em tautological} holomorphic line bundle over the base of the family.  The Hermitian connection and Chern form for the renormalized metric are determined.  Connections to the work of M. Mirzakhani, L. Takhtajan and P. Zograf, and intersection numbers for the moduli space of punctured Riemann surfaces studied by E. Witten are presented.
\end{abstract}

\section{Comparing cusps}
The renormalization of a hyperbolic metric at a cusp is introduced.  The setting is used to present an intrinsic norm for the germ of a holomorphic map at a cusp.  

A compact Riemann surface $R$ having punctures and negative Euler characteristic has a complete hyperbolic metric, \cite{Ahconf}. The geometry of a cusp of a hyperbolic metric is standard.  From the uniformization theorem for a puncture $p$ there is a distinguished local conformal coordinate with $z(p)=0$ and the metric locally given by the germ of
\begin{equation}
\label{cusp}
ds^2=\left(\frac{|dz|}{|z|\log |z|}\right)^2.
\end{equation}
The distinguished coordinate is given as $z=e^{2\pi i\zeta}$ for the cusp represented at infinity for the upper half plane $\mathbb H$  with coordinate $\zeta$ and translation $\zeta\rightarrow \zeta+1$ .  
The {\em canonical} coordinate $z$ is unique modulo a unimodular factor and the circle $|z|=c$ is the closed horocycle about $p$ with hyperbolic length $\ell=- 2\pi/\log c$.  The complete hyperbolic metric for the punctured unit disc $\mathbb D=\{0<|z|<1\}$ is likewise given by formula (\ref{cusp}).  The unit area neighborhood of the puncture $\{\ell(z)\le 1\}\subset\mathbb D$ isometrically embeds to a neighborhood of each puncture of a complete hyperbolic metric, \cite{Leut}.

A conformal coordinate $\xi$ for a neighborhood of a puncture $p$ can be used to describe the conformal completion at $p$.  The {\em omitted point} $p$ corresponds to the {\em omitted value} of $\xi$.  The tangent space at $p$ is described in terms of the tangent space at the omitted value of 
$\xi$.  We use these general observations to define a norm for the tangent space at the puncture.  
\begin{definition}
Let $R$ be a Riemann surface with hyperbolic metric and a puncture $p$ with canonical local coordinate $z$ as above.  The canonical norm for the tangent space at $p$ is defined by $\|\frac{\partial}{\partial z}\|_{can}=1$.
\end{definition}

As noted the coordinate $z$ is unique modulo a unimodular factor.  For canonical coordinates $z$ and $w$ respectively for complete hyperbolic metrics $ds_1^2$ and $ds_2^2$ consider a germ at the origin of a holomorphic map $w=h(z)=zg(z)$ with $h(0)=0$. We have for the ratio of metrics
\begin{multline}
\label{comparison}
(h^*ds_2^2)(ds_1^2)^{-1}=\left(\frac{|zh'(z)|\log|z|}{|h(z)|\log|h(z)|}\right)^2 = \left(\left|\frac{h'(z)}{g(z)}\right|\frac{\log |z|}{\log |zg(z)|}\right)^2\\
=1\,-\,2\,\frac{\log|h'(0)|}{\log|z|}\,+\,O_g\left(\frac{1}{(\log |z|)^2}\right).
\end{multline}

We use the observations to present a formula for the norm of a map between tangent spaces.  Consider $R_1$ and $R_2$ with complete hyperbolic metrics $ds_1^2$ and $ds_2^2$ with punctures $p_1$ and $p_2$.  For the germ $h$ of a holomorphic map from a neighborhood of $p_1$ to a neighborhood of $p_2$ the canonical norms for the tangent spaces determine a norm for the differential $dh$.  Observe that a point $q$ close to $p_1$ lies on a unique simple closed horocycle about the puncture.  For $q$ close to $p_1$ write $\ell(q)$ for the associated horocycle length.  
\begin{lemma}
\label{norm}
Let $h$ be the germ of a holomorphic map between punctures.  At the puncture the differential satisfies $\lim\limits_{q\rightarrow p_1}(h^*ds_2^2)(ds_1^2)^{-1}=1$ and $\log\|dh\|_{can}^2=\lim\limits_{q\rightarrow p_1} \frac{2\pi}{\ell}\log ((h^*ds_2^2)(ds_1^2)^{-1})$.  For the map $h$ expressed in terms of canonical coordinates $w=h(z)$ the norm satisfies $\|dh\|_{can}=|h'(0)|.$
\end{lemma}
{\bf Proof.} For canonical coordinates at the punctures the formulas follow from the observation $\ell=-2\pi/\log |z|$ and expansion (\ref{comparison}).  The proof is complete.

We will study the variation of $\lim \frac{2\pi}{\ell}\log(h^*ds_2^2)(ds_1^2)^{-1}$ for a family of Riemann surfaces.  In the following sections we consider holomorphic families of Riemann surfaces.  In preview, for $\nu(s)$ a suitable family of Beltrami differentials holomorphic in a parameter $s$, $\nu(s)$ with fixed compact support, a holomorphic family $\{R^{\nu(s)}\}$ is defined.  A conformal coordinate $\zeta$  for a neighborhood of a puncture $p$ of $R$ also serves as a conformal coordinate for $R^{\nu(s)}$ provided the support of the coordinate is disjoint from $supp(\nu(s))$. Accordingly the tangent vector $\frac{\partial}{\partial\zeta}$ is a holomorphic section of the family of tangent spaces along $p$.  In this case the hyperbolic metrics are locally conformal and easily compared as $ds_{R^{\nu(s)}}^2=e^{2f}ds_R^2$.  For the canonical coordinates at $p$, $z$ for $R$ and $w$ for $R^{\nu(s)}$, the composition of  $w(\zeta)\circ(z(\zeta)^{-1})$ is the germ of a holomorphic map $h$. We have from Lemma \ref{comparison} that 
$\log \|dh\|_{can}^2=\lim\limits_{q\rightarrow p}\frac{4\pi}{\ell}f$.  We will study the limit in the following sections. 

We separately note that the Schwarz lemma also provides a comparison for germs of hyperbolic metrics at a puncture, \cite{Ahconf}.  For $\mathcal U$ a punctured neighborhood of a puncture $p$ with germs of hyperbolic metrics $ds_1^2$ and $ds_2^2$, let $\mathcal V \subset \mathcal U$ be a neighborhood with $\partial \mathcal V$ compact in $\mathcal U$.  The ratio $\rho=ds_2^2/ds_1^2$ has limit unity at $p$ from the comparison (\ref{comparison}).  In a local conformal coordinate $\zeta$ the constant curvature equation is $-(ds^2)^{-1}\Delta\log ds^2=-1$ for $ds^2$ now the metric local expression and $\Delta$ the coordinate Euclidean Laplacian.  The pair of metrics satisfy $\Delta \log ds^2_2/ds_1^2=ds_2^2-ds^2_1$.  The ratio $\rho$ is continuous and has a maximum on $\bar{\mathcal V}\cup\{p\}$.  Provided the maximum of $\rho$ is greater than unity, the difference of metrics is positive at the maximum.  In this case from the combined equations $\log \rho$ is strictly convex at the maximum and the maximum cannot be interior.   It follows that the maximum of $\rho$ is either unity or occurs on $\partial \mathcal V$.  An  application is provided for $\mathcal V$ the $ds_1^2$ horoball neighborhoods $\{\ell(q)<\ell\}$ of the cusp.  The maximum of $\rho$ on $\{\ell(q)\le \ell\}$ is either unity or occurs only on $\{\ell(q)=\ell\}$. A general property follows.  The maximum $\max_{\{\ell(q)=\ell\}}\rho$  either has the constant value unity for all small $\ell$, or is strictly increasing in $\ell$. The analogous statement for $\min_{\{\ell(q)=\ell\}}\rho$ is also valid.

\section{The prescribed curvature solution}
The solution of the prescribed curvature equation describes the ratio of a pair of hyperbolic metrics in a neighborhood of a puncture. For a family of hyperbolic metrics we consider the variation in parameters of the prescribed curvature solution.  For calculation of the connection and Chern forms it is required that variations of the prescribed curvature solutions are suitably bounded at the cusps.  We find that adapted variations are suitably bounded. Our considerations begin with the basics of variations of conformal structures and variations of solutions of the prescribed curvature equation.  We provide the required formulas and estimates.     

Let $R$ be a Riemann surface with hyperbolic metric $ds^2$ and canonical line bundle $\kappa$.  Let $S(p,q)$, 
$p,q\in \mathbb Z$ be the space of smooth sections of $\kappa^{p/2}\otimes\overline{\kappa^{q/2}}$ 
and $S(r)=S(r,-r)$.  For a $p$-differential $\varphi\in S(2p,0)$ the product $(ds^2)^{-p/2}\varphi$ 
is an element of $S(p)$ and the conjugate $(ds^2)^{-p/2}\overline{\varphi}$ an element of $S(-p)$.  An element of 
$S(r)$ has a well-defined absolute value.  Metric-derivative operators (essentially the covariant derivatives) on $S(r)$ are defined as follows: for $\zeta$ a local conformal coordinate and $ds^2=\rho^2(\zeta)|d\zeta|^2$ then $K_r=\rho^{r-1}\frac{\partial}{\partial\zeta}\rho^{-r}$ and 
$L_r=\rho^{-r-1}\frac{\partial}{\partial\overline{\zeta}}\rho^{r}$, \cite{Wlhyp}.  Note for $f\in S(r)$ that $\overline{f}\in S(-r)$, $K_r:S(r)\rightarrow S(r+1)$, $L_r:S(r)\rightarrow S(r-1)$ and $K_r=\overline{L_{-r}}$.  
Further note for $f\in S(p)$, $g\in S(r-p)$ that $K_r(fg)=gK_pf+fK_{r-p}g$ and $L_r(fg)=gL_pf+fL_{r-p}g$.  The hyperbolic metric Laplacian operator on functions is given as $D=4L_1K_0$. For $f\in S(r)$ the absolute value $|f|$ is a function and a $C^k$ norm is defined  by $\|f\|_0=\sup |f|$ and $\|f\|_k=\sum_{|P|\le k}\|Pf\|_0$ for the sum over all products of operators $K_*$ and $L_*$ of length at most $k$.  We are also interested in functions suitably small in the cusps.  A weighted-normed space $C_{\ell^n}$ is defined for $f\in S(r)$ by $\|f\|_{\ell^n}=\sup_{\ell\le 1}|\ell^{-n}f|\vee\sup_{R-\{\ell\le 1\}}|f|$.  In effect from Lemma \ref{comparison} we are interested in the variational formulas for a family of hyperbolic metrics with variational terms bounded in $C_{\ell}$. 

Let $R$ be a punctured Riemann surface and $Q(R)$ the associated space of holomorphic quadratic differentials with at most simple poles at punctures.  The space $B(R)$ of continuous Beltrami differentials is the finite $C^0$-norm elements of $S(-2)$.  A basic feature of Riemann surfaces is the integral pairing $\int_R \mu\varphi$ for $\mu\in B(r)$ and $\varphi\in Q(R)$.  The deformation space of $R$ is the Teichm\"{u}ller space $\mathcal T$, \cite{Ahsome, Earlet, Ng}. We will write $\mathcal T_{g,n}$ for the Teichm\"{u}ller space of genus $g$, $n$ punctured surfaces.  At $R$ the holomorphic tangent space to the deformation space is $B(R)\slash B(R)^{\perp}$ and the holomorphic cotangent space is $Q(R)$.  

A new conformal structure $R^{\mu}$ for $\mu\in B(R)$ with $\|\mu\|<1$ is defined as follows.  For $z$ a local conformal coordinate for $R$ and $w(z)$ a local homeomorphism solution on the range of $z$ for the Beltrami equation $w_{\bar z}=\mu w_z$ the composition $w\circ z$ is a local conformal coordinate for $R^{\mu}$.  In brief for $\{z_{\alpha}\}_{\alpha \in \mathcal A}$ an atlas for the $R$-conformal structure, the compositions $\{w_{\alpha}\circ z_{\alpha}\}_{\alpha \in \mathcal A}$ provide an atlas for the $R^{\mu}$-conformal structure.  Local holomorphic coordinates for $\mathcal T$ are described as follows: for $\{\nu_1,\dots,\nu_m\}$ continuous Beltrami differentials spanning $B(R)\slash B(R)^{\perp}$ and $s\in \mathbb C^m$, set $\nu(s)=\sum_j s_j\nu_j$; for $s$ small $s \rightarrow R^{\nu(s)}$ is a local holomorphic coordinate for $\mathcal T$, \cite{Ahsome, Earlet,Ng}.

The harmonic Beltrami differentials $\mathcal H(R)=\{\mu\in B(R)\mid K_{-2}\mu =0\}$ form a natural subspace.  The interest with harmonic Beltrami differentials arises from the observation that variational formulas are simplified by their introduction, \cite{Ahsome,Ng, Wfharm, Wlchern}.  Basic properties are as follows: $\mathcal H(R)\subset B(R)$ is a direct summand for $B(R)^{\perp}$ and consequently $\mathcal H(R)\simeq B(R)\slash B(R)^{\perp}$; the mapping $\varphi \rightarrow(ds^2)^{-1}\bar\varphi$ is a complex anti linear bijection of $Q(R)$ to $\mathcal H(R)$.

The considerations of the first section require maps holomorphic in a neighborhood of the punctures.  To this purpose we introduce smooth truncations of elements of $\mathcal H(R)$.  Begin with a choice of local coordinate $z$ for a neighborhood $\mathcal U$ of the puncture $p$ with $z(p)=0$.  Choose a function $\chi$ on $\mathbb C$ to be an approximate characteristic function of $\{|z|>1\}$ such that $\chi$ is smooth, vanishes identically on $\{|z|\le 1/2\}$ and is identically unity on $\{|z|\ge 1\}$.  For $\epsilon$ sufficiently small the $\epsilon$-{\em truncation} of $\mu\in \mathcal H(R)$ at the puncture $p$ is locally defined as $\chi(z/\epsilon)\mu(z)$; as $\epsilon$ tends to zero the truncation tends to $\mu$.  The truncation of $\mu$ is defined by introducing a truncation at each puncture.  For $\mu$ harmonic $K_{-2}\chi\mu=\mu K_0\chi$ and we combine local expressions at a puncture for the hyperbolic metric and a holomorphic quadratic differential to find that $\|K_{-2}\chi\mu\|_1$ is $O(\epsilon \,\log^4 \epsilon)$ for $\epsilon$ tending to zero. 

To compare metrics, we introduce a particular {\em pullback} of a metric.  For a Riemann surface $S$ with almost complex structure $J_S$, $J_S^2=-id$; arbitrary metric $d\sigma^2$ and K\"{a}hler form $\omega$ we note that $d\sigma^2=2\omega(\ ,J_S\ )$.  For a Riemann surface $R$ with almost complex structure $J_R$, $J_R^2=-id$, and $h:R\rightarrow S$ a smooth map the $J$-{\em pullback by} $h$  metric is defined as $d\sigma_*^2=2h^*\omega(\ , J_R\ )$.  The $J$-pullback metric has a straightforward local coordinate description.  For $z$ a conformal coordinate for $R$, $w$ a conformal coordinate for $S$ and $d\sigma^2=(\alpha(w)|dw|)^2$ then $d\sigma_*^2=\alpha(w(z))^2(|w_z|^2-|w_{\bar z}|^2)|dz|^2$.  As explained in \cite[pg. 449]{Wlhyp}, introduction of the $J$-pullback metric simplifies various variational formulas and in particular gives rise to formulas which are {\em local} in the Beltrami differential and its derivatives.  By comparison introduction of the standard pullback $h^*d\sigma^2$ entails solving the potential equation $F_{\bar z}=\mu$.

To compare hyperbolic metrics we begin with $ds_R^2$ the hyperbolic metric of $R$ and a smooth map $h:R\rightarrow S$. Introduce now $d\sigma_*^2$ the $J$-pullback by $h^{-1}$ metric of $ds_R^2$ to 
$S$. (Use of the inverse map will simplify the resulting expressions.)  The metric $d\sigma_*^2$ on $S$ is a perturbation of the hyperbolic metric on $R$.  Let $C$ be the Gaussian curvature of $d\sigma_*^2$ and $C_*=C\circ h$.  Further let $D_{\sigma_*}$ be the $d\sigma_*^2$-Laplacian on $S$ and $D_*$ its pullback by $h$, i.e. for $g$ a smooth function on $S$ then 
$(D_{\sigma_*}g)\circ h=D_*(g\circ h)$.  

We review the formulas for $C_*$ and $D_*$  needed for the present considerations. The complete variational formulas  were presented in \cite[Sec. 5.2, esp. Lemma 5.3]{Wlhyp}.   We will consider a family  $\{S=R^{\mu(a)}\}$ for $\mu(a)=a\mu$ for a Beltrami differential $\mu\in B(R)$ and a  real parameter $a$.
The curvature $C_*(a)$ and pullback-Laplacian $D_*(a)$ depend real analytically on $a$, \cite{Wlhyp}.  Basic properties are as follows.  In general for $supp(\mu)\subset \mathcal C$, we have on the complement $R-\mathcal C$ that $C_*(a)=-1$ and (for local expressions) $D_*(a)=D$.    The variational derivatives of $C_*$ and $D_*$ for $\mu(a)=a\mu$ at $a=0$ are as follows
\begin{equation}
\label{CandD1}
C_a=4\Re K_{-1}K_{-2}\mu, \quad C_{aa}=-D|\mu|^2-8\Re(L_1(\bar\mu K_{-2}\mu))
\end{equation}
and for $\mu$ denoting the multiplication operator 
\begin{equation}
\label{CandD2}
D_a=-8\Re K_{-1}\mu K_0, \qquad \mbox{\cite[Example 5.4]{Wlhyp}}.
\end{equation}
The formulas simplify for harmonic Beltrami differentials since $K_{-2}\mu$ vanishes and are similarly approximated for truncated differentials since $K_{-2}\mu$ is $O(\epsilon \,\log^4 \epsilon)$.

We are ready for $h:R\rightarrow S$ to compare metrics, beginning with $ds_R^2$ the complete hyperbolic 
metric of $R$.  The complete hyperbolic metric of $S$ is  described in terms of the $J$-pullback by $h^{-1}$ metric $d\sigma_*^2$ and the {\em solution of the prescribed curvature equation}
\[
ds_S^2=e^{2\tilde f}d\sigma_*^2 \quad \mbox{for} \quad D_{\sigma_*}\tilde f-C=e^{2\tilde f},\quad \mbox{\cite[pg. 453]{Wlhyp}}.
\]
As noted above the composition $f=\tilde f\circ h$ satisfies the equation $D_*f-C_*=e^{2f}$.  The equation $D_{\sigma_*}\tilde f-C=e^{2\tilde f}$ is equivalent to the curvature of $e^{2\tilde f}d\sigma_*^2$ being identically $-1$, \cite{Wlhyp}. The uniformization theorem and estimates of the first section provide that the prescribed curvature equation has a unique solution in $C_{\ell}$, \cite{Ahconf}.  Standard results from Teichm\"{u}ller theory provide that for $\mu(a)\in B(R)$ varying real analytically the prescribed curvature solution $f(a)$ varies real analytically in $C^k$-norm for compact subsets of $R$, \cite{AB,Earlet}.  

We are now ready to consider the variation of the prescribed curvature solution.  The expressions for the Laplacian $D_*(a)$, curvature $C_*(a)$ and solution $f(a)$ are real analytic in the parameter $a$.  The expressions for the variations in $a$ are determined.  We have {\em for small values of} $a$ the equations
\begin{equation}
\label{fa}
D_*f-e^{2f}=C_*, \qquad   
D_*f_a+D_af-2e^{2f}f_a=C_a
\end{equation}
and after adding $Df_a$ to each side of the second equation
\begin{equation}
\label{fa*}
Df_a-2e^{2f}f_a=(D-D_*)f_a-D_af+C_a.
\end{equation}
For a Beltrami differential $\mu$ with compact support the right hand side of the third equation  has compact support.    We now use the equations to show that $f(a)$ is a differentiable $C_{\ell}$-valued function.  For $\ell$ the horocycle length function, the $C_{\ell^n}$ weighted-norm is $\|f\|_{\ell^n}=\sup_{\ell\le 1}|\ell^{-n}f|\vee\sup_{R-\{\ell\le 1\}}|f|$.
\begin{lemma}
\label{Cell}
Let $f(a)$ be the prescribed curvature solution for hyperbolic metrics and Beltrami differential $\mu(a)=a\mu$, $\mu$ with compact support. The solution $f(a)$ is a differentiable $C_{\ell}$-valued function of the parameter $a$.  For small values of $a$ the solution satisfies
\[
f(a)=af_a(0)+O(a^2)\quad and \quad f_a(a)=f_a(0)+af_{aa}(0)+O(a^2)
\]
for remainder terms bounded in $C_{\ell}$. 
\end{lemma}
{\bf Proof.} We iterate basic bounds to establish the desired expansions.  In preview, a bound for $f$ gives a bound for $e^{2f}$ which from equation (\ref{fa*}) gives a bound for the derivative $f_a$ and from integration in the parameter an improved bound for $f$.  
We start by noting from the final paragraph of the first section that $\|f\|_0$ is bounded by $\sup_{supp(\mu)}|f|$ and as noted above on $supp(\mu)$ the solution $f$ is real analytic in $C^k$ in $a$. Lemmas 11 and 12 below will be used to pass from estimates on $supp(\mu)$ for the right hand side of (\ref{fa*}) to estimates on R. To begin, from Lemma \ref{D-g inverse} and the bounds for quantities on the compact support of the right hand side of (\ref{fa*}) it follows from Lemma \ref{Cl2toCl1} that $f_a$ is uniformly bounded in $C_{\ell}$ for all suitably small $a$.  We next integrate  in the parameter $\int_{a'}^{a''}f_a\,da$ to find for all small $a$, $a-a'$ that $f(a)$ is $O(a)$ in $C_{\ell}$ and that $f(a)=f(a')+O(|a-a'|)$ in $C_{\ell}$. The expansion is the Taylor expansion for a continuous $C_{\ell}$-valued function. 

We next compare the equations (\ref{fa*}) for values $a$ and $a'$ with $e^{2f(a)}=e^{2f(a')+O(|a-a'|)}$ and apply Lemmas \ref{D-g inverse} and \ref{Cl2toCl1} to find for all small $a$, $a-a'$ that $f_a(a)=f_a(a')+O(|a-a'|)$ in $C_{\ell}$.  We next integrate in the parameter to find that $f(a)=f(a')+(a-a')f_a(a')+O((a-a')^2)$ in $C_{\ell}$. The expansion is the Taylor expansion for a differentiable $C_{\ell}$-valued function.  

We next substitute the first-order expansions for the right hand side $rhs(a)$ of (\ref{fa*}) and also the expansion $e^{2f(a)}=1 + 2af_a(0) +O(a^2)$ in $C_{\ell}$ to find the equation in $C_{\ell}$ (modulo terms of order $a^2$ in $C_{\ell^2}$) for $f_a(a)$
\[
Df_a-2f_a-2af_af_a=rhs(0)+a\,rhs_a(0)+ O_{C_{\ell^2}}(a^2).
\]
The resulting equation
\[
f_a=(D-2)^{-1}\left(rhs(0)+a\,rhs_a(0)+2af_a(0)^2+O_{C_{\ell^2}}(a^2)\right)
\]
provides for the desired Taylor expansion for $f_a(a)$. We note in passing that all remainder terms depend on $supp(\mu)$.  The proof is complete.

\section{The canonical curvature}

We present the formula for the Hermitian connection and curvature of the canonical norm for the family of tangent spaces along a puncture.  The curvature is given in terms of the Takhtajan-Zograf form.

Let $\Gamma$ be a Fuchsian group uniformizing the Riemann surface $R$ with $n$ punctures.  Let $\Gamma_1,\dots,\Gamma_n$ be a set of non-conjugate maximal parabolic subgroups of $\Gamma$ representing the punctures.  Further let $\sigma_j\in PSL(2;\mathbb R)$ be transformations such that $\sigma_j^{-1}\Gamma_j\sigma_j=\Gamma_{\infty}$ with $\Gamma_{\infty}$ the group of integer-translations.  The Eisenstein series for the $j^{th}$ puncture of $R$ is defined for $\zeta\in\mathbb H$ and $\Re s>1$ as
\[
E_j(\zeta,s)=\sum_{\gamma\in\Gamma_j\backslash\Gamma}\Im (\sigma_j^{-1}\gamma\zeta)^s,\qquad j=1,\dots,n,\ \cite{Boraut, Vn}.
\]
The special value $E_j(\zeta)=E_j(\zeta,2)$ plays an important role for deformation theory, \cite{WlO,TZ,Wldis}.  

The Fourier expansions of the Eisenstein series have a simple form, \cite{Boraut, Vn}.  Conjugate the $k^{th}$ cusp to infinity and 
consider $E_j(\sigma_k(\zeta))=E_j(\sigma_k(\zeta),2)$ which is $\Gamma_{\infty}$ invariant.  For $\zeta=x+iy$, the Kronecker delta $\delta$, and coefficients $c_{jk}$, there is the expansion for large $y$
\[
E_j(\sigma_k(\zeta))=\delta_{jk}y^2+c_{jk}y^{-1}+O(e^{-2\pi y}).
\]
The functions $E_j(\zeta)$ satisfy the differential equation $DE=2E$ and are not elements of $L^2(R)$.  

The hyperbolic metric Laplacian $D$ acting on $L^2(R)$ has non positive spectrum.  The Green's function $G(\zeta,\zeta')$ for the operator $(D-2)^{-1}$ is a part of our considerations. The Green's function is the integral kernel with respect to the hyperbolic area element $dA$ for inverting the $(D-2)$ operator: for $u(\zeta)=\int_RG(\zeta,\zeta')g(\zeta')dA$ then $(D-2)u=g$ for $g\in C^0$.  We again conjugate the $j^{th}$ cusp to infinity and consider Fourier expansions.  For $\zeta'$ fixed, $\zeta=x+i y$, there is the expansion for $y$ large
\[
G(\sigma_j(\zeta),\zeta')=c_j(\zeta')y^{-1}+O(e^{-2\pi y})
\]
and for $g$ with compact support
\[
u(\sigma_j(\zeta))=c_j(u)y^{-1}+O(e^{-2\pi y}).
\]

We begin considerations with the approach of L. Takhtajan and P. Zograf to express the leading coefficient of $(D-2)^{-1}g$ in terms of the integral $\int_R gE\,dA$, \cite[Lemma 2]{TZ}.  The approach follows the argument for the Maass-Selberg relation, \cite{Boraut, Vn}.
\begin{lemma}
\label{limit}
Let $\Gamma$ be a Fuchsian group uniformizing a Riemann surface with punctures.  For $g$ with compact support and $\ell_j$ the horocycle length function for the $j^{th}$ puncture, the operator $(D-2)^{-1}$ satisfies \[
\lim_{\Im \zeta\rightarrow\infty}\ell_j^{-1}((D-2)^{-1}g)(\sigma_j\zeta)=\frac{-1}{3}\int_R gE_j\,dA.
\]
\end{lemma}
{\bf Proof.}  We begin with the defining property $\hat g=(D-2)\hat u$ and observe that 
\[\int_RE_j\circ\sigma_j g\, dA=\int_R E_j\circ\sigma_j Du-2E_j\circ\sigma_j u\,dA=\int_RE_j\circ\sigma_j Du-DE_j\circ\sigma_j u\,dA
\]
(where we write $g$ for $\hat g\circ\sigma_j$ and $u$ for $\hat u\circ\sigma_j$.)  The terms of the third integrand are not individually integrable since the leading coefficient of $E_j\circ\sigma_j u$ at infinity is $\Im \zeta$.  To consider the third integral we introduce a fundamental domain $\mathcal F$ for $\sigma_j^{-1}\Gamma\sigma_j$, containing the cusp neighborhood $\{0\le \Re \zeta <1, \Im \zeta>1\}$.  We further introduce the  sub domains $\mathcal F^Y=\{\zeta\in\mathcal F\mid \Im\zeta\le Y\}$ and from the above
\[
\int_RE_j\circ\sigma_jg\,dA=\lim_{\Im Y\rightarrow\infty}\int_{\mathcal F^Y}E_j\circ\sigma_jDu-DE_j\circ\sigma_ju\,dA
\]
and apply Green's formula to find
\[
=\lim_{\Im Y\rightarrow\infty}\int_{\partial \mathcal F^Y}E_j\circ\sigma_j\frac{\partial u}{\partial\mathbf n}-
\frac{\partial E_j\circ\sigma_j}{\partial\mathbf n}u\,ds
\]
for the hyperbolic metric elements $\frac{\partial}{\partial\mathbf n}$ and $ds$.   The integral over 
$\partial\mathcal F^Y \cap\partial\mathcal F$ vanishes by consideration of the orientation and the group invariance of $E_j\circ\sigma_j$ and $u$.  The Fourier expansions for $E_j\circ\sigma_j$ and $u$ provide that the remaining integral over $\{0\le\Re \zeta <1,\Im \zeta=Y \}$ is $-3\,c_j(\hat u) + O(e^{- Y})$. The proof is complete. 

The tangent space to the deformation space at a Riemann surface $R$ is represented by the space $\mathcal H(R)$ of harmonic Beltrami differentials.  A Hermitian form for $\mathcal H(R)$ defines a Hermitian metric for Teichm\"{u}ller space.  For $\mu,\nu\in\mathcal H(R)$ and hyperbolic area element $dA$ the Weil-Petersson (WP) form is
\[
\langle\mu,\nu\rangle_{WP}=\int_R\mu\bar\nu\, dA, \qquad \cite{Ahsome}
\]
and for $R$ with punctures $p_1,\dots,p_n$ the Takhtajan-Zograf (TZ) form for the puncture $p_j$ is
\[
\langle\mu,\nu\rangle_{TZ,p_j}=\int_R\mu\bar\nu E_j\,dA, \qquad \cite{TZ}.
\]
The Takhtajan-Zograf metric is $\sum_j\langle\ ,\ \rangle_{TZ,p_j}$.  We are ready to present the formula for the Hermitian connection and curvature for the canonical norm.
\begin{theorem}
\label{main}
Let $\mathcal T$ be the Teichm\"{u}ller space of a Riemann surface $R$ with punctures.  The canonical norm $\|\ \|_{can,p}$ for the family of tangent spaces along the puncture $p$ has Hermitian connection vanishing on $\mathcal H(R)$ and Chern form $c_1(\|\ \|_{can,p})=\frac{2i}{3}\overline{\langle\ ,\ \rangle}_{TZ,p}$ on $\mathcal H(R)$.
\end{theorem}
{\bf Proof.} We combine considerations.  We have from the discussion following 
Lemma \ref{norm} that for the family $\{R^{\nu(s)}\}$, $\nu(s)$ a truncated Beltrami differential with compact support, and $\zeta$ a local conformal coordinate at a puncture, that  
$\log \|\frac{\partial}{\partial \zeta}\|^2_{can}=\lim_{q\rightarrow p}\frac{4\pi}{\ell}f$ for $ds^2_{R^{\nu(s)}}=e^{2f}ds_R^2$.  The hyperbolic metrics of $\{R^{\nu(s)}\}$ are described in terms of the solutions $f(s)$ of the prescribed curvature equation.  The Taylor expansions of $f(s)$ of Lemma \ref{Cell} provide that $\lim_{q\rightarrow p}\frac{4\pi}{\ell}f(s)$ is $C^2$ at $s=0$ with initial $s$-derivative $\lim_{q\rightarrow p}\frac{4\pi}{\ell}f_s(0)$ and initial $s\bar s$-derivative $\lim_{q\rightarrow p}\frac{4\pi}{\ell}f_{s\bar s}(0)$.   From equation (\ref{fa}) and Lemma \ref {limit} we find in terms of the curvature, Laplacian and prescribed curvature solution the variations at $s=0$
\[
\lim_{q\rightarrow p}\frac{4\pi}{\ell}f_s(0)=\lim_{q\rightarrow p}\frac{4\pi}{\ell}(D-2)^{-1}C_s=\frac{-4\pi}{3}\int_RC_sE_p\,dA
\]
and
\begin{multline*}
\lim_{q\rightarrow p}\frac{4\pi}{\ell}f_{s\bar s}(0)=
\lim_{q\rightarrow p}\frac{4\pi}{\ell}(D-2)^{-1}(C_{s\bar s}-D_sf_{\bar s}-D_{\bar s}f_s+4f_sf_{\bar s}) \\
=\frac{-4\pi}{3}\int_R (C_{s\bar s}-D_sf_{\bar s}-D_{\bar s}f_s+4f_sf_{\bar s})E_p\,dA
\end{multline*}
(the formulas for $f_s(0)$ and $f_{s\bar s}(0)$ are a straightforward calculation.)  We next observe that since the canonical norm is twice differentiable its Hermitian connection $1$-form $\Theta$ and curvature $2$-form $\Omega$ have well-defined evaluations on $\mathcal H(R)$.  The evaluations are given by the limit of the truncated Beltrami differentials introduced in the prior section (the limit as the parameter $\epsilon$ tends to zero.)  We have from the estimate for truncated  Beltrami differentials combined with (\ref{CandD1}) and (\ref{CandD2}) for elements of $\mathcal H(R)$ the resulting formulas at $s=0$: $C_s=0$, $f_s=(D-2)^{-1}C_s=0$ and $C_{s\bar s}=\frac{-1}{2}D|\mu|^2$.  As the final step we apply Green's formula to note that 
$\frac{-1}{2}\int_RD|\mu|^2E_p\,dA=-\int_R|\mu|^2E_p\,dA$ and finally that: 
\[
\Theta(\mu)=0,\ \Omega(\bar\mu,\mu)=\frac{4\pi}{3}\overline{\langle\mu,\mu\rangle}_{TZ,p}\ 
\mbox{ and }\ c_1(\|\ \|_{can,p})=\frac{i}{2\pi}\Omega=\frac{2i}{3}\overline{\langle\ ,\ \rangle}_{TZ,p}.
\]
The proof is complete.

Associated to a Hermitian form $\langle\ ,\ \rangle$ is a (pre) K\"{a}hler form $\frac{i}{2}\langle\ ,\ \rangle$. The above result provides a new proof that the TZ metric is K\"{a}hler, \cite{TZ}.  We now write 
$\omega_{TZ,p}=\frac{i}{2}\langle\ ,\ \rangle$ and restate the result.
\begin{corollary}
\label{tgtcurv}
The Chern form of the canonical norm for the family of tangent spaces along a puncture $p$ satisfies $c_1(\|\ \|_{can,p})=\frac{-4}{3}\omega_{TZ,p}$.
\end{corollary}
Takhtajan-Zograf were able to determine the Chern form for a single puncture without identifying the metric \cite[formula (11)]{TZ}, while L. Weng using Arakelov theory determined the Chern form for multiple punctures without identifying the metric \cite{Weng}.

\section{Applications}

We present beginning properties of the canonical norm Chern form and connections to the work of other authors.  The TZ metric is K\"{a}hler, \cite{TZ} and incomplete, \cite{Obit1}.  K. Obitsu, W. K. To and L. Weng \cite{OTW} have recently determined the asymptotic behavior of the metric akin to the original result of H. Masur, \cite{Msext}.  A simple property of the metric comes from the observation 
$E=\sum_{\gamma\in\Gamma_{\infty}\backslash\Gamma}\Im (\gamma z)^2$ and that the integral $\int_R\mu\bar\nu E\,dA$ can be {\em unfolded} to $\int_{\Gamma_{\infty}\backslash\mathbb H}\mu\bar\nu(\Im \zeta)^2\,dA$ which is a special value of the Rankin-Selberg convolution $L$-function, \cite[Theorem 2]{TZ}.

The families of cotangent spaces along punctures give rise to rational cohomology classes on the moduli space $\mathcal M_{g,n}$ of genus $g$, $n$ punctured Riemann surfaces.  By definition $\psi_j$ is the rational Chern class of the orbifold line bundle whose fiber at the point $[R;p_1,\dots,p_n]\in \mathcal M_{g,n}$ is the cotangent space at $p_j$ (see \cite{ArCor} for the definition on $\overline{\mathcal M_{g,n}}$.)  The {\em canonical norm} provides a metric for the line bundles $\psi_j$; for the canonical local coordinate $z$ at a puncture, $\|dz\|_{can}=1$.  We restate our main formula in the present setting.  
\begin{corollary}
\label{chern}
The Chern form of the canonical norm for the family of cotangent spaces $\psi_p$ along a puncture $p$ satisfies $c_1(\|\ \|_{can,p})=\frac43\omega_{TZ,p}$.
\end{corollary}
L. Takhtajan and P. Zograf used Quillen's metric to calculate the first Chern form of the determinant line bundle for families of $\bar\partial$-operators, \cite{TZ}.  The authors considered the Teichm\"{u}ller space $\mathcal T_{g,n}$, the Teichm\"{u}ller curve $\mathcal C_{g,n}\rightarrow\mathcal T_{g,n}$ and $\mathcal E_k=T^{-k}_{vert}\,\mathcal C_{g,n}$ the $k^{th}$ symmetric power of the dual of the vertical line bundle of $\mathcal C_{g,n}\rightarrow\mathcal T_{g,n}$.  On a fiber $\mathcal C_{g,n}\rightarrow\mathcal T_{g,n}$ of the Teichm\"{u}ller curve the vertical tangent space $T_{vert}\,\mathcal C_{g,n}$ coincides with the tangent space of the fiber.  Associated to the family of $\bar\partial$-operators for $\mathcal E_k$ is an index bundle $\ind \bar\partial_k$ and determinant holomorphic line bundle 
$\det\ind\bar\partial_k$ with a Quillen metric $\|\ \|_{Quillen}$ determined from the hyperbolic metric of Riemann surfaces.   Takhtajan-Zograf found a local index formula using \cite{Wlchern} for families of compact Riemann surfaces with $2g-2>0$, $k\ge 0$,
\[
c_1(\det\ind\bar\partial_k)=\frac{6k^2-6k+1}{12\pi^2}\omega_{WP}, \quad \mbox{\cite{TZcomp}}
\]
and for families of punctured Riemann surfaces with $2g-2+n>0$, $k\ge 0$,
\[
c_1(\det\ind\bar\partial_k)=\frac{6k^2-6k+1}{12\pi^2}\omega_{WP}-\frac19\sum_j\omega_{TZ,p_j},\quad\mbox{\cite{TZ}.}
\]
Quillen's metric involves the zeta function determinant of the Laplacian.  For punctured Riemann surfaces 
Takhtajan-Zograf used special values of the Selberg zeta function in place of zeta function determinants, \cite[formula (6)]{TZ}.  

In \cite{Wlhyp} the family hyperbolic metric for the vertical line bundle  $T_{vert}\,\mathcal C_g\rightarrow\mathcal T_g$ was used to find a Chern form on $\mathcal C_g$ and to calculate the pushdown of the square of the form.  We found that the pushdown class $\kappa_1$ is represented by the pushdown form $\frac{1}{\pi^2}\omega_{WP}$ (see the section below on the WP K\"{a}hler form.)  By using truncated harmonic Beltrami differentials the formula can be generalized to families of punctured Riemann surfaces.   We now combine results and present a local form of the above Takhtajan-Zograf formula.   
\begin{corollary}
For bundles over $\mathcal T_{g,n}$  the Quillen metric, vertical line bundle metric and cotangent spaces along punctures metric determined from the hyperbolic metric there is a pointwise relation of Chern forms
\[
12\,c_1(\det\ind\bar\partial_k)=(6k^2-6k+1)\,c_1(\kappa_1)-\sum_jc_1(\psi_{p_j}).
\]
\end{corollary}
Certain comments are in order.  First, the hyperbolic metric is determined by a choice of conformal structure and does not involve a choice of {\em marking} and so the above considerations are valid for $\mathcal M_{g,n}$ the moduli space of punctured Riemann surfaces.  Second, the considerations for the compactified moduli space of stable curves $\overline{\mathcal M_{g,n}}$ have not been effected. 

An application is the curvature of the conormal bundle to the divisor of noded Riemann surfaces.  Families of cotangent spaces along punctures can be used to describe the conormal bundle.  In particular a pair of families of punctured Riemann surfaces $\{R\}$ and $\{R'\}$ and a formal pairing of the punctures $p$ of $R$ and $p'$ of $R'$ determines a family $\mathcal D=\{R\vee R'\}$ of noded Riemann surfaces, \cite{Bersdeg}.  Consider the family $\mathcal M$ where the node $p\vee p'$ is allowed to open.  The product of cotangent spaces along $p$ and $p'$ defines a line bundle $\lambda$ over $\mathcal D\subset \mathcal M$.  We recall that $\lambda$ is isomorphic to the conormal bundle of the divisor $\mathcal D\subset\mathcal M$.  The family $\mathcal M$ of noded Riemann surfaces can be described in terms of deformations supported away from the node and the {\em plumbing family} $\{(z,w,t)\mid zw=t,\,|z|,|w|,|t|<1\}\rightarrow\{|t|<1\}$.  For such a description  the function $t$ becomes a local defining function for the divisor $\mathcal D\subset \mathcal M$.  
For a change of parameterization $f(z),g(w)$ and $h(t)$ for the plumbing family  with $f(0),g(0)$ and $h(0)$ each zero there is the basic relation $f'(0)g'(0)=h'(0)$.  The relation provides the cocycle relation for the isomorphism of the conormal bundle $\lambda$ and the product of cotangent spaces along punctures.  The product of canonical norms for cotangent spaces provides a norm for $\lambda$, as well as for the inverse bundle $\lambda^{-1}$.  From Corollary 
\ref{tgtcurv} the curvature $c_1(\lambda^{-1})$ is negative definite and we find a local form of the principle 
that {\em the opening of a node is negative}.  

L. Weng studied for punctured Riemann surfaces the {\em intersection product} for metrized line bundles and also the Deligne-Riemann-Roch isometry, \cite{Weng}.  He introduced the metrized WP, TZ and logarithmic Mumford line bundles over $\mathcal M_{g,n}$ and determined first Chern forms.  As part of his results he showed for the metrized TZ line bundle $\underline{\Delta_{TZ}}$ on $\mathcal M_{g,n}$ that $c_1(\underline{\Delta_{TZ}})=\frac43\omega_{TZ}$, \cite[pg. 278]{Weng}.

An application is for the volume of moduli spaces.  M. Mirzakhani has considered the moduli space $\mathcal M_{g,n}(b_1,\dots,b_n)$ of genus $g$ bordered Riemann surfaces with geodesic boundary components of prescribed length $(b_1,\dots, b_n)$, \cite{Mirvol}.  The K\"{a}hler form $\omega_{WP}=\frac 12 \sum d\ell\wedge d\tau$ provides a symplectic form on $\mathcal M_{g,n}(b_1,\dots,b_n)$.  Mirzakhani developed a recursive scheme for determining the volumes.  Using an identity for geodesic length she established a general volume-result, \cite{Mirvol}.  
\begin{theorem}
The volume $V_{g,n}(b)=Vol(\mathcal M_{g,n}(b_1,\dots,b_n))$ is a polynomial in the squares of geodesic boundary lengths $b_1^2,\dots,b_n^2$ with 
\[
V_{g,n}(b)=\sum_{|\alpha|\le 3g-3+n}c_g(\alpha)b^{2\alpha}
\]
where $\alpha$ ranges over multi indices of $(\mathbb Z_{\ge 0})^n$ and $c_g(\alpha)$ are positive values of $\pi^{6g-6+2n-2|\alpha|}\mathbb Q$.
\end{theorem}
Mirzakhani applied the result in \cite{Mirvol} to provide a volume-expansion for tubular neighborhoods of the compactification divisor $\mathcal D$ in the Deligne-Mumford compactified moduli space $\overline{\mathcal M}$.  

The orbifolds $\mathcal M_{g,n}(b_1,\dots,b_n)$ form an $(\mathbb R_{\ge 0})^n$ bundle over $\mathcal M_{g,n}$ the moduli space of genus $g$, $n$ punctured Riemann surfaces.  A Riemann surface with geodesic boundaries  and a {\em point} on each boundary is alternately described by an $n$ punctured Riemann surface and a product of $n$ factors of $S^1$, a principal torus-bundle over a punctured Riemann surface.  Mirzakhani finds in \cite{Mirwitt} that symplectic reduction (for an $(S^1)^n$ quasi-free action following Guillemin-Sternberg \cite{Gubk}) can be used to provide a simple description for the family of K\"{a}hler forms
\begin{equation}
\label{symred}
\omega_{WP}\big\vert_{\mathcal M_{g,n}(b_1,\dots,b_n)}=\omega_{WP}\big\vert_{\mathcal M_{g,n}}+ \sum\limits_j \frac{b_j^2}{4}\,c_1(\psi_j)
\end{equation}
(the relation is for cohomology classes on $\overline{\mathcal M}$.)  (Mirzakhani considers the symplectic form $2\,\omega_{WP}$ and so formulas differ by a factor of $2$.)  The expansion is presented for a general choice of principal connection for the $S^1$ bundles.  An explicit principal connection is given by introducing the line bundles $\psi_j$, the canonical norms $\|\ \|_{can,j}$ and Hermitian connections.  The consequence is a local form of the above expansion.  The local expansion agrees with a perturbation formula of K. Obitsu and the author, \cite{WlO}. We considered the perturbation of the $WP$ metric and K\"{a}hler form for the tangent subspaces parallel to the compactification divisor $\mathcal D\subset \overline{\mathcal M}$.    The expansion is a refinement to the work of H. Masur \cite{Msext}, G. Daskalopoulos and R. Wentworth \cite{DW2}, and the author \cite{Wlcomp}.  For a family $\{R_{\ell}\}$ of hyperbolic surfaces given by {\em pinching} short geodesics all with common length $\ell$, we found for the K\"{a}hler forms restricted to the tangent subspaces parallel to the compactification divisor
\[
\omega_{WP}^{tgt}(\ell)= \omega_{WP}^{tgt}(0)+ \frac{\ell^2}{3}\sum\limits_j\omega_{TZ,p_j}(0)+O(\ell^3)
\]
and with Corollary \ref{chern} the pointwise relation
\[
= \omega_{WP}^{tgt}(0) + \frac{\ell^2}{4}\sum\limits_j\,c_1(\psi_{p_j})+O(\ell^3).
\]

Mirzakhani combined her integration scheme and formula (\ref{symred}) to show that the collection of integrals (intersection pairings)
\[
\int_{\overline{\mathcal M_{g,n}}}c_1(\psi_1)^{\alpha_1}\cdots c_1(\psi_n)^{\alpha_n}\omega_{WP}^{3g-3+n-|\alpha|}
\]
satisfies the recursion for the {\em string equation} and the {\em dilaton equation}. The intersection numbers combine to provide a {\em partition function} $F$ for two-dimensional quantum gravity, \cite{Witt1,Witt2}.  E. Witten conjectured that $e^F$  would satisfy the $KdV$ equations $\mathbf L_k\,e^F=0$, $k\ge-1$ with Virasoro constraint relations $[\mathbf L_m,\mathbf L_k]=(m-k)\mathbf L_{m+k}$.  M. E. Kazarian and S. K. Lando \cite{KaLa}, Y.-S. Kim and K. Liu \cite{KiLi}, M. Kontsevich \cite{Kont}, M. Mulase and B. Safnuk \cite{MuSa}, A. Okounkov and R. Pandharipande \cite{OP}, and Mirzakhani \cite{Mirwitt} have verified the conjecture.  The authors show that the Virasoro relations determine the intersection numbers of tautological line bundles.  Mirzakhani's integration scheme also determines the intersection numbers.  A consequence of Corollary \ref{chern} is the following. 
\begin{corollary}
The Virasoro relations or Mirzakhani's integration scheme can be used to determine all TZ-WP pairings 
\[\int_{\overline{\mathcal M_{g,n}}}\omega_{TZ,1}^{\alpha_1}\cdots\omega_{TZ,n}^{\alpha_n}\,\omega_{WP}^{3g-3+n-|\alpha|}.
\]
\end{corollary}

We illustrate the result with an example.  Mirzakhani provided in \cite{Mirvol} the expansion for $2\,\omega_{WP}$ for $(g,n)=(0,4)$
\[
V_{0,4}(b)=\frac12(4\pi^2+b_1^2+\cdots+b_4^2)
\]
which corresponds to the integrals
\[
\pi^2=\int_{\overline{\mathcal M_{0,4}}}\omega_{WP}=\pi^2\int_{\overline{\mathcal M_{0,4}}}\kappa_1
\]
using that $\kappa_1=\frac{1}{\pi^2}\omega_{WP}$ and from (\ref{symred}) to the integrals
\[
\int_{\overline{\mathcal M_{0,4}}}c_1(\|\ \|_{can,j})=1.
\]
The values agree with the evaluations of E. Arbarello and M. Cornalba \cite{ArCor} and P. Zograf \cite{Zogsph}.  We apply the last evaluation to find the TZ volume of $\overline{\mathcal M_{0,4}}$.  The volume form is $dV_{TZ}=\sum_j\omega_{TZ,j}$ and consequently
\[
\int_{\overline{\mathcal M_{0,4}}}dV_{TZ}=3.  
\]


\section{The WP K\"{a}hler form}

The characteristic class of the WP K\"{a}hler form is part of the present considerations.  We now revisit our earlier treatment of the WP symplectic and K\"{a}hler forms to find that 
$\omega_{WP, symplectic}=2\,\omega_{WP, K\ddot{a}hler}$.  Certain earlier formulas especially for integrals and characteristic classes need to be adjusted for the present considerations.

An underlying real tangent space $V$ for a complex manifold $M$ has an almost complex structure $J$, $J^2=-id$.  The complexification $V^{\mathbb C}$ of the tangent space $V$ is decomposed into the $\pm i$ eigenspaces of $J$.  The decomposition is given as $V^{1,0}\oplus V^{0,1}$ with $V\subset V^{\mathbb C}$ the subspace fixed by complex conjugation. For a Riemann surface $R$ at the corresponding point of the Teichm\"{u}ller space  the holomorphic tangent space is $V^{1,0}\simeq \mathcal H(R)$ and the holomorphic cotangent space is $(V^{1,0})^*\simeq Q(R)$.  The  tangent-cotangent pairing for $\mu$ in $\mathcal H(R)$ and $\varphi$ in $Q(R)$ is $(\mu,\varphi)=\int_R \mu\varphi$.  On the holomorphic tangent space $V^{1,0}\simeq\mathcal H(R)$ the WP Hermitian form is $\langle\ ,\ \rangle_{WP}$ and the K\"{a}hler form $\omega_{WP}=\frac{i}{2}\langle\ ,\ \rangle_{WP}$ (corresponding to the K\"{a}hler form $\frac{i}{2}(dz_1\wedge d\bar z_1\cdots dz_m\wedge d \bar z_m)$ on $\mathbb C^m$.) We will see below that in effect $i\langle\ ,\ \rangle_{WP}$ was used in our earlier papers \cite{WlFN, Wlsymp, Wlpi, Wlhomol, Wlchern, Wlhyp}.  Takhtajan-Zograf use the K\"{a}hler form $\frac{i}{2}\langle\ ,\ \rangle_{WP}$, \cite[pg. 402]{TZ}.

We studied in \cite{WlFN} for a closed geodesic $\alpha$ on a Riemann surface the relationship between the Fenchel-Nielsen infinitesimal twist deformation $t_{\alpha}$ and the geodesic-length function $\ell_{\alpha}$. For $\theta_{\alpha}$ the classical Petersson theta series for $\alpha$, the infinitesimal Fenchel-Nielsen twist $t_{\alpha}$ is represented by the harmonic Beltrami differential $\frac{i}{\pi}(dA)^{-1}\overline{\theta_{\alpha}}$ in $\mathcal H(R)$ ($dA$ the hyperbolic area element), \cite[Corollary 2.8]{WlFN}.   F. Gardiner's formula for the differential of geodesic-length is $d\ell_{\alpha}=\frac{2}{\pi}\Re\int_R \mu\theta_{\alpha}$, \cite{Gardtheta}.   A K\"{a}hler form $\frac{i}{2}\langle\ ,\ \rangle$ is evaluated with a sum over permutation of vectors. We calculate the twist-length duality
\[
(\mu, d\ell_{\alpha})=\frac{2}{\pi}\Re\int_R \mu\theta_{\alpha}=2\Re i\langle\mu,\frac{i}{\pi}(dA)^{-1}\overline{\theta_{\alpha}}\rangle_{WP}=2\,\omega_{WP,K\ddot{a}hler}(\mu,t_{\alpha}).
\]
The twist-length duality in terms of real tangent vectors is $d\ell_{\alpha}=2\,\omega_{WP,K\ddot{a}hler}(\ ,t_{\alpha})$.  If we write $\omega_{WP,symplectic}$ for the symplectic form used in our earlier papers then we have $\omega_{WP,symplectic}=2\,\omega_{WP,K\ddot{a}hler}(\ ,\ )$, \cite[see Theorem 2.10]{WlFN}. The definition of the symplectic form and twist-length duality formula continued in our subsequent papers.  There are corresponding adjustments to subsequent formulas.  Integration formulas are relevant for the present considerations.  The $\mathcal M_{1,1}$ and $\mathcal M_{0,4}$ area formulas become $\int_{\mathcal M_{1,1}}2\,\omega_{WP,K\ddot{a}hler}=\frac{\pi^2}{6}$ and 
$\int_{\mathcal M_{0,4}}2\,\omega_{WP,K\ddot{a}hler}=2\pi^2$ since the twist-length duality formula was used to derive the integrand.  To represent the appropriate characteristic class there is a further factor of $\frac12$ for $\mathcal M_{1,1}$ since the WP pairing should be for integration over tori modulo their elliptic involution.  The {\em canonical coordinates} formula becomes 
$\omega_{WP,K\ddot{a}hler}=\frac12\sum_j d\ell_j\wedge d\tau_j $ \cite[Theorem 1.3]{Wldtau}.  The characteristic class formula \cite[formula (5.1)]{Wlhomol} should also be adjusted for the definition of the K\"{a}hler form.  The formula for $\kappa_1$ the pushdown of the square of the Chern form for the family hyperbolic metric is likewise affected \cite[see proof of Corollary 5.11]{Wlchern} and \cite{Wlhyp}.  The updated formula $\kappa_1=\frac{1}{\pi^2}\omega_{WP,K\ddot{a}hler}$ for \cite{Wlchern} and updated formula (5.1) for \cite{Wlhomol} now combine to agree with the formula $\kappa_1=12\lambda -\delta $ of D. Mumford, 
\cite{Mumlens} and the  Takhtajan-Zograf $\kappa_1$ calculation \cite[pg. 424]{TZ}.

\section{Estimates for Green's operators}
We provide estimates for the operator $(D-k)^{-1}$ acting on functions small at the cusps.

\begin{lemma}
\label{D-g inverse}
For $k\in C^0$ with a positive infimum $m_k$  the operator $(-D+k)^{-1}$ is continuous on $L^2(R)$.  For a continuous $g\in L^2(R)$ the operator satisfies $|(-D+k)^{-1}g|\le(-D+m_k)^{-1}|g|$.
\end{lemma}
{\bf Proof.}  For $\beta$ the supremum of $k$ we write $k=\beta-\hat k$ and introduce the factorization $(-D+k)^{-1}=(1-(-D+\beta)^{-1}\hat k)^{-1}(-D+\beta)^{-1}$.  By hypothesis $m_k=\beta - \sup \hat k$ is positive and we can consider the geometric series
\begin{equation}
\label{g series}
(1-(-D+\beta)^{-1}\hat k)^{-1}= 1 + (-D+\beta)^{-1}\hat k+((-D+\beta)^{-1}\hat k)^2+\cdots 
\end{equation}
($\hat k$ is now the multiplication operator.)  From the spectral theorem the $L^2$-norm of the operator $(-D+\beta)^{-1}$ is $\beta^{-1}$ and the above series of $L^2$-operators converges since $\beta^{-1}\sup\hat k < 1$. The operator $(-D+k)^{-1}$ is defined on $L^2$.  Next we consider the pointwise behavior of $(-D+c)^{-1}g$, for a 
continuous $g\in L^2$.  The integral kernel for $(-D+c)^{-1}$, $c$ a positive constant, is positive and so for a non negative continuous function $v\in L^2$ we have the inequality $(-D+\beta)^{-1}\hat k\,v\le (-D+\beta)^{-1}(\beta- m_k)\,v$.  To estimate $(-D+k)^{-1}g$ we apply the right hand side of (\ref{g series}) to  $(-D+\beta)^{-1}|g|\ge |(-D+\beta)^{-1}g|$ and apply the above inequality for each factor of $(-D+\beta)^{-1}\hat k$.  The result is a convergent geometric series expansion for $(-D+m_k)^{-1}|g|$, the desired upper bound.  The proof is complete. 

We now consider the operator $(D-2)^{-1}$ acting on functions in $C_{\ell^n},\,n>1$.  The Green's function for the operator $(-D+2)$ is given as an absolutely convergent sum
\[
G(z,z_0)=\sum_{\gamma\in\Gamma}Q(z,\gamma z_0)
\]
for $\delta(z,z_0)$ hyperbolic distance and $Q(z,z_0)=-Q_2(\delta(z,z_0))$ for $Q_2$ an associated Legendre function, \cite[Chap. 1]{Fay}.  We consider that $\Gamma$ contains the group of integer translations $\Gamma_{\infty}$ as a maximal parabolic subgroup and introduce
\[
G_{\infty}(z,z_0)=\sum_{\gamma\in\Gamma_{\infty}}Q(z,\gamma z_0) \quad\mbox{and}\quad 
G_{\dagger}(z,z_0)=G(z,z_0)-G_{\infty}(z,z_0).
\]
We write $\mathcal C=\{\ell\le 1\}\subset R$ for the union of the unit area horocycle regions and $\mathcal C_{\infty}$ for the unit area horocycle region at infinity.  In the following we analyze the behavior of $\int_R G(z,z_0)g(z_0)\,dA$ by considering the region of integration as the union of $R-\mathcal C$, $\mathcal C-\mathcal C_{\infty}$ and finally $\mathcal C_{\infty}$.  The initial estimate was presented in \cite{Wlhyp}.  In Lemma A.4.1 of \cite[pg. 468]{Wlhyp} we showed given $\delta_0>0$ there exists a constant $c_0$ such that $0<G(z,z_0)\le c_0\, e^{-\delta(z,z_0)}$ for hyperbolic distance $\delta(z,z_0)>1$ provided the injectivity radius at $z$ or $z_0$ is at least $\delta_0$.  

\begin{lemma}
\label{Cl2toCl1} 
For $n>1$, the operator $(D-2)^{-1}$ is a continuous mapping from $C_{\ell^n}$ to $C_{\ell}$.
\end{lemma}
{\bf Proof.} We begin with the contribution to $\int Gg\,dA$ for the compact region $R-\mathcal C$.  Lemma A.4.1, \cite{Wlhyp}, provides that the integral is valued in $C_{\ell}$ since for $z_0\in\mathcal C$, $z\in R-\mathcal C$  then $e^{-\delta(z,z_0)}\le \ell$. The bound for the compact region is complete. 

The estimates for the cusp regions are indicated by considering the  region $\mathcal C_{\infty}$ at infinity.  We consider for $z\in\mathcal C_{\infty}$ the contribution
\[
\int_{\mathcal C-\mathcal C_{\infty}}G(z,z_0)g(z_0)\,dA\ +\ \int_{\mathcal C_{\infty}}G_{\dagger}(z,z_0)g(z_0)\,dA.
\]
For the first integral $z\in\mathcal C_{\infty}$ and $z_0\in\mathcal C-\mathcal C_{\infty}$, while in the second integral $z,z_0\in \mathcal C_{\infty}$.  On the indicated regions $G$ and $G_{\dagger}$ are smooth solutions of the differential equation $(D-2)u=0$.  For $z$ fixed and $z_0$ tending to a puncture from Lemma A.4.1 each function tends to zero.  It follows from the maximum principle that for $z$ fixed $G$ and $G_{\dagger}$ achieve their maxima in $z_0$ respectively on $\partial(\mathcal C-\mathcal C_{\infty})$ and $\partial\mathcal C_{\infty}$.  The estimate of Lemma A.4.1 applied for $z_0$ respectively on $\partial(\mathcal C-\mathcal C_{\infty})$ and $\partial\mathcal C_{\infty}$ now provides a uniform bound for $G$ and $G_{\dagger}$ by a constant multiple of the horocycle length function $\ell$.  Since $g$ is integrable on $R$ the desired bounds for the integrals follow. 

It remains to consider the principal contribution $\int_{\mathcal C_{\infty}}G_{\infty}(z,z_0)g(z_0)\,dA$. For the upper half plane $\mathbb H$ the integral is {\em unfolded} and $g$ is replaced by $y^{-n}$ to give the integral
\begin{equation}
\label{integral}
\int_{\Im z_0\ge 1}Q(z,z_0)y^{-n-2}dxdy
\end{equation}
for the variable $z_0=x+i y$.  We introduce a majorant for $Q$.  The values of the function $y^{-1}$ on a disc $\{\delta(z,z_0)\le 1\}$ are within fixed multiples of the value at the center.  It follows for sake of bounding (\ref{integral}) that for hyperbolic distance $\delta\le 1$ the contribution of the kernel $Q$ is bounded by a comparison kernel with a positive minimum on the disc.  It now follows from the large hyperbolic distance description of the behavior  $-Q_2(\delta(z,z_0))\le ce^{-2\delta(z,z_0)}$, \cite[pg. 468]{Wlhyp}, and the formula $\cosh \delta(z,z_0)=1+\frac{|z-z_0|^2}{2\Im z\Im z_0}$ that the integral (\ref{integral}) is bounded for $a=\Im z$ in terms of

\begin{multline*}
\int_{\substack{-\infty<x<\infty \\ y\ge 1}}\left(1+\frac{|z-ia|^2}{2ya}\right)^{-2}y^{-n-2}dxdy \\
=4\int_{\substack{-\infty<x<\infty \\ y\ge 1}}\frac{a^2}{y^n(x^2+y^2+a^2)^2}dxdy=2\pi\int_{y\ge 1}\frac{a^2}{y^n(y^2+a^2)^{\frac32}}dy.
\end{multline*}
Upon substituting $\frac{a^2}{(y^2+a^2)^{\frac32}}\le \frac{1}{a}$ the last integral is bounded by 
$O(a^{-1})$. The integral (\ref{integral}) is bounded by a multiple of the horocycle length function $\ell=a^{-1}$, as desired.  The proof is complete.


\end{document}